\documentclass[12pt]{article}
\usepackage{amssymb, amsmath}
\begin{document}
\voffset=0.0truein \hoffset=-0.5truein \setlength{\textwidth}{6.0in}
\setlength{\textheight}{8.8in} \setlength{\topmargin}{-0.2in}
\renewcommand{\theequation}{\arabic{section}.\arabic{equation}}
\newtheorem{thm}{Theorem}[section]
\newtheorem{defin}{Definition}[section]
\newtheorem{lem}{Lemma}[section]
\newtheorem{prop}{Proposition}[section]
\newtheorem{cor}{Corollary}[section]
\newcommand{\n}{\nonumber}
\newcommand{\w}{\omega}
\newcommand{\s}{\sigma}
\newcommand{\bn}{\|_{\dot{B}^{\frac{N}{p}}_{p,1}}}
\renewcommand{\b}{\dot{B}^{\frac{N}{p}+1}_{p,1}}
\newcommand{\ib}{B^{\frac{N}{p}+1}_{p,1}}
\renewcommand{\th}{\theta}
\newcommand{\la}{\lambda}
\newcommand{\La}{\Lambda}
\newcommand{\g}{\gamma}
\renewcommand{\O}{\Omega}
\newcommand{\tv}{\tilde{v}}
\newcommand{\tw}{\tilde{\omega}}
\renewcommand{\a}{\alpha}
\renewcommand{\o}{\omega}
\newcommand{\e}{\varepsilon}
\renewcommand{\t}{\theta}
\newcommand{\nao}{\nabla ^\bot \theta}
\newcommand{\vare}{\varepsilon}
\newcommand{\bb}{\begin{equation}}
\newcommand{\ee}{\end{equation}}
\newcommand{\bq}{\begin{eqnarray}}
\newcommand{\eq}{\end{eqnarray}}
\newcommand{\bqn}{\begin{eqnarray*}}
\newcommand{\eqn}{\end{eqnarray*}}
\title{On the behaviors of solution near possible blow-up time in the incompressible Euler and related equations }
\author{Dongho Chae\thanks{ This research is done, while the author was visiting University of Chicago.
The work  is supported partially by KRF Grant(MOEHRD, Basic Research Promotion Fund).}\\
Department of Mathematics\\
              Sungkyunkwan University\\
               Suwon 440-746, Korea\\
              {\it e-mail : chae@skku.edu}}
 \date{}
\maketitle
\begin{abstract}
We study behaviors of scalar quantities near the possible blow-up time, which is made of smooth solutions of the Euler equations, Navier-Stokes equations and the surface quasi-geostrophic equations. Integrating the dynamical equations of the scaling invariant norms,  we derive the possible blow-up behaviors of the above quantities, from which we obtain new type of blow-up criteria and some necessary conditions for the blow-up.
\end{abstract}

\section{ The  Main Theorems}
 \setcounter{equation}{0}
\subsection{The Euler equations}
 We are  mainly concerned with the following Euler equations for the
homogeneous incompressible fluid flows in $\Bbb R^N,$  $N\geq 2$.
$$
(E)\left\{\aligned
&\partial_t v +(v\cdot \nabla )v=-\nabla \pi ,\\
&\textrm{div }\, v =0,\\
&v(x,0)=v_0(x).
\endaligned
\right.
$$
Here $v=(v_1, \cdots, v_N )$, $v_j =v_j (x, t)$, $j=1, \cdots, N$, is the
velocity of the flow, $\pi=\pi(x,t)$ is the scalar pressure, and $v_0 $
is the given initial velocity, satisfying div $v_0 =0$.
Given $k\in \Bbb N\cup \{0\}$ and $p\in [1, \infty)$  we use $W^{k,p} (\Bbb R^N)$ to denote the standard Sobolev space on $\Bbb R^N$,
$$
W^{k,p} (\Bbb R^N)=\left\{ f\in L^p (\Bbb R^N)\,\Big|\, \sum_{j=0}^k \|D^j f\|_{L^p} <\infty\right\},
$$
and set $W^{k,2} (\Bbb R^N)=H^k (\Bbb R^N).$
We also use notations for the solenoidal vector fields in Sobolev space,
$$W^{k,p}_\sigma (\Bbb R^N)=\left\{ v=(v_1, \cdots, v_N)\,\Big|\, v_j\in W^{k,p} (\Bbb R^N), j=1,\cdots ,N, \,\, \mathrm{div}\, v =0\right\},
$$
and $H^k_\sigma (\Bbb R^N)=W^{k,2}_\sigma (\Bbb R^N),$ $L^p_\sigma (\Bbb R^N)=W^{0,p}_\sigma (\Bbb R^N)$.
Given $v_0\in H^k _\s(\Bbb R^N)$, $k>N/2 +1$, there exists local in time unique solution $v\in C([0, T);H^k_\s (\Bbb R^N))$ for some $T=T(\|v_0 \|_{H^k})$(\cite{kat1}). The question that if this local smooth solution can be continued arbitrary longer time or not  is one of the most outstanding open problem in the mathematical fluid mechanics(see e.g. \cite{maj, che} for an introduction to the subject, and \cite{con5, bar,cha5} for more recent survey article). If the local smooth solution cannot be continued beyond $T_* <\infty$, then $\lim\sup_{t\to T_*} \|v(t)\|_{H^m} =\infty$, in which case we say the blow-up happens at finite time $T_*$. As one direction of research for the  problem  people have been trying to derive sharp blow-up criterion, which was initiated by Beale-Kato-Majda(\cite{bea}), which shows that the integral $\int_0 ^{T_*} \|\o (t)\|_{L^\infty} dt$(vorticity magnitude) controls blow-up at $t=T_*$, where $\o=$ curl $v$ is the vorticity(see e.g. \cite{koz1,koz2,cha0} for later refinements in this direction).  Later Constantin-Fefferman-Majda(\cite{con3, con1}, which was developed from \cite{con2})  took into account geometric structure of the vortex stretching term in the vorticity equations to get another kind of blow-up condition, where the dynamics of the direction of vorticity play essential roles(see also \cite{den1,den2} for later refinements in this direction)  These two separate forms of criteria controlling the blow-up by magnitude and the direction of the vorticity respectively are interpolated in \cite{cha1}, developing the observations in \cite{cha2}. In the current paper we investigate the behaviors of smooth solutions near the possible blow-up time to derive different type of below-up criteria from previous ones.  For this purpose we introduce suitable scalar quantities made of smooth solutions of the equations.
The main feature of equations we use is the scaling invariance properties of the Euler system, and the analysis of self-similar form of equations is essential for our result. In this aspect it is a development of the author's previous studies of the possibility of self-similar blow-up  of the Euler equations(\cite{cha3,cha4}). We also derive  blow-up conditions in the Navier-Stokes and the surface quasi-geostrophic equations in the next subsections.\\
Let $\xi (x,t)=\o(x,t)/|\o(x,t)|$ be the direction field of the vorticity for 3D vector field $v(x,t)$, and $S=(S_{jk})$ with $S_{jk}=\frac12 (\partial_jv_k +\partial_k v_j )$ be the deformation tensor. In \cite{con3} the following useful (local) quantity was defined.
\bb\label{loc}
\a (x,t)=\sum_{j,k=1}^3\xi _j (x,t) S_{jk}(x,t)\xi_k(x,t) .
\ee
In the case $\o(x,t)=0$ we set $\a(x,t)=0$.
For a nonzero smooth solution $v(x,t)\in C([0, T); H^k _\s(\Bbb R^N))$ with $k\in \Bbb N \cup\{0\} $, we define
the following (nonlocal) quantity
\bb\label{glo}
\a_{k}(t) :=\frac{-\int_{\Bbb R^N} D^k \left[(v(x,t)\cdot \nabla)v(x,t) \right] \cdot D^k v(x,t) dx}{\|D^k v (t)\|_{L^2}^2}.
\ee
If $v(\cdot, t)=0$, then we set $\a_k (t)=0$.
In the following we derive necessary and sufficient conditions for blow-up.
\begin{thm}
Let  $k> N/2+1$, and  $v\in C([0, T_*); H^k _\s(\Bbb R^N))$ be a local smooth solution to (E) with $v_0\in H^k _\s(\Bbb R^N)$, $v_0\neq 0$. Then, the following (i)-(iii) are necessary and sufficient conditions for blow-up at $T_*$.
\begin{itemize}
\item[(i)] There exists an absolute constant $K=K(N,k)$ such that
\bb\label{1.2}
\lim\inf_{t\to T_*} (T_*-t)\|D^k v(t)\|_{L^2}^{\frac{N+2}{2k}} \geq \frac{K}{\|v_0\|_{L^2}^{1-\frac{N+2}{2k}}}.
\ee
\item[(ii)]
\bb\label{1.3}
\lim\inf_{t\to T_*} \int_{0}^{t} \left\{\a_k (\tau)-\left[\frac{2k}{N+2}\right]\frac{1}{T_* -\tau}\right\}d\tau >-\infty.
\ee
\item[(iii)] For all $\vare_0 >1$ there exists a sequence $\{t_n \}$ with $t_n \nearrow T_*$ such that
\bb\label{1.4}
\a_k (t_n)\geq \left[\frac{2k}{N+2}\right]\frac{1}{T_* -t_n}-\left[\frac{2k\vare_0}{N+2}\right]\frac{1 }{T_* -t_n} \left[\log \left(\frac{1}{T_* -t_n}\right)\right]^{-1}
\ee
for all $n=1,2,\cdots.$
\item[(iv)]$(N=3)$ For all $\vare_0 >1$ there exists a sequence $\{t_n \}$ with $t_n \nearrow T_*$ such that
\bb\label{1.4a}
\|\a (t_n)\|_{L^\infty}\geq \frac{1}{T_* -t_n}-\frac{\vare_0 }{T_* -t_n} \left[\log \left(\frac{1}{T_* -t_n}\right)\right]^{-1}
\ee
for all $n=1,2,\cdots.$
\end{itemize}
\end{thm}
{\it Remark 1.1 } Since
\bqn
\lefteqn{\left|\int_{\Bbb R^N} D^k \left[(v(x,t)\cdot \nabla)v(x,t) \right] \cdot D^k v(x,t) dx\right|}\hspace{.4in}\n \\
&&=\left|\int_{\Bbb R^N}\left\{ D^k \left[(v(x,t)\cdot \nabla)v(x,t) \right]- (v\cdot \nabla )D^k v \right\}\cdot D^k v(x,t) dx\right|\n \\
&&\leq \int_{\Bbb R^N} |D^k \left[(v(x,t)\cdot \nabla)v(x,t) \right]- (v\cdot \nabla )D^k v | |D^k v(x,t)| dx\n \\
&&\leq \|D^k \left[(v(x,t)\cdot \nabla)v(x,t) \right]- (v\cdot \nabla )D^k v \|_{L^2} \|D^k v\|_{L^2}\n \\
&&\leq \hat{C}\|\nabla v\|_{L^\infty} \|D^k v \|_{L^2} ^2
\eqn
(see (\ref{2.37}) and the first part of the proof of Theorem 1.1 below), we have
\bb\label{1.5}
|\a_k (t)|\leq \hat{C}\|\nabla v\|_{L^\infty}
\ee
for an absolute constant $\hat{C}=\hat{C}(k, N)$.
Hence, (\ref{1.4}) implies that a necessary and sufficient condition of blow-up at $T_*$ is that there exists a constant $C$ depending on $k, N$ such that
\bb\label{1.6}
\lim\sup_{t\to T_*}(T_*-t)\|\nabla v(t)\|_{L^\infty} \geq \hat{C}\lim\sup_{t\to T_*}(T_*-t)\a_k (t)\geq \frac{2k\hat{C}}{N+2}.
\ee
In the case $N=3$, using the obvious inequality
\bb\|\a (t)\|_{L^\infty} \leq \|\nabla v (t)\|_{L^\infty}
\ee
instead of (\ref{1.5}), then  similarly to the above we can deduce from (\ref{1.4a})
\bb\label{1.6a}
\lim\sup_{t\to T_*}(T_*-t)\|\nabla v(t)\|_{L^\infty} \geq 1 ,
\ee
which is also obtained in \cite{cha6} in a different context.\\
\ \\
In order to state the next theorem we recall that a positive continuous function $g(\cdot)$ is said to satisfy  {\em Osgood's condition}
if
$$ \int_1 ^\infty \frac{ds}{g(s)} <\infty.
$$
The Osgood condition is a necessary and sufficient condition for finite time blow-up of the ordinary differential equation with positive initial data(\cite{osg}),
$$ \frac{dy}{dt}= g(y,t), \qquad y(0)=y_0 >0.
$$

\begin{thm} Let $k> N/2+1$, and $v\in C([0, T_*); H^k_\s (\Bbb R^N))$ be a local smooth solution to (E) with $v_0\in H^k_\s (\Bbb R^N)$, $v_0\neq 0$.
If $T_*$ is the  blow-up time, then  one of the following three statements hold true.
\begin{itemize}
\item[(i)]  There exists a sequence $\{ t_n \}_{1=0}^\infty$ with $t_n \nearrow T_*$ such that
\bb\label{1.7}
\a _{k} (t_n)=\left[\frac{2k}{N+2}\right]\frac{1}{T_* -t_n} \qquad \forall n=1,2,\cdots .
\ee

\item[(ii)] There exists $t_0 \in (0, T_*)$ such that
\bb
 \left|\a_k (t)-\left[\frac{2k}{N+2}\right] \frac{1}{T_*-t}\right| >0 \qquad \forall t\in (t_0, T_*),
\ee
and
\bb \label{1.9}
\lim_{t\to T_*} (T_*-t)\|D^k v(t) \|_{L^2}^{\frac{N+2}{2k}}+\int_{0}^{T_*} \left|\a_{k}(t)-\left[\frac{2k}{N+2}\right]\frac{1}{T_*-t}\right|dt <\infty.
  \ee
  \item[(iii)]
  There exists $t_0 \in (0, T_*)$ such that
\bb\label{1.10}
\a_k (t)>\left[\frac{2k}{N+2}\right] \frac{1}{T_*-t}\qquad \forall t\in (t_0, T_*),
\ee
  and
 \bb\label{1.11}
\lim_{t\to T_*}(T_*-t)\|D^k v(t) \|_{L^2}^{\frac{N+2}{2k}}=\int_{0}^{T_*} \left\{\a_{k}(t)-\left[\frac{2k}{N+2}\right]\frac{1}{T_*-t}\right\}dt=\infty.
\ee
Furthermore, for any continuous, positive function $g$ satisfying the Osgood condition, we have
 \bb\label{1.12}
 \int_{t_1} ^{T_*} \frac{\left\{ \a_{k}(t)-\left[\frac{2k}{N+2}\right]\frac{1}{T_* -t}\right\}}{g\left(\log\{(T_*-t)\|D^k v(t) \|_{L^2}^{\frac{N+2}{2k}}\|v_0\|_{L^2}^{1-\frac{N+2}{2k}} \}\right)} dt <\infty
\ee
for  $t_1\in (0, T_*)$ sufficiently close to $T_*$.
\end{itemize}
\end{thm}
{\it Remark 1.2 } We note that the integrability condition in (\ref{1.9}) implies that there exists a sequence $\{ t_n \}$ with $t_n \nearrow T_*$ such that
\bb\label{1.13}
\a _{k} (t_n)=\left[\frac{2k}{N+2}\right]\frac{1}{T_* -t_n} + o \left(\frac{1}{T_* -t_n }\left[\log \left(\frac{1}{T_* -t_n }\right)\right]^{-1}\right)  \quad\mbox{as $n \to \infty$},
\ee
which is a special case of (\ref{1.4}).
Indeed, we can write the integral in (\ref{1.9}) as
 \bb
  \int_{0}^{T_*} \left\{\log \left( \frac{1}{T_* -t}\right)\left|(T_*-t)\a_{k} (t)-\frac{2k}{N+2}\right|\right\}\frac{1}{(T_* -t)\log\left( \frac{1}{T_* -t}\right) }dt <\infty.
  \ee
  Since
  $$ \int_{0} ^{T_*} \frac{1}{(T_* -t)\log \left( \frac{1}{T_* -t}\right) }dt =\infty, $$
  there should be a sequence $\{ t_n \}$ with $t_n \nearrow T_*$ such that
  \bb
  \lim_{n\to \infty}\log \left( \frac{1}{T_* -t_n}\right)\left|(T_*-t_n)\a_{k} (t_n)-\frac{2k}{N+2}\right|=0.
  \ee
 We observe that (\ref{1.7}) belongs to the case of (\ref{1.13}). We can thus further narrow down the possibilities for the behavior of $\a_k(t)$ near the possible blow-up time $T_*$ as follows:
 \\ Either (\ref{1.13}) holds for a sequence $t_n \nearrow T_*$, or (\ref{1.10}) holds
for some $t_0<T_*$. This is more specified than Theorem 1.1(iii).\\

\subsection{The Navier-Stokes equations}
In this subsection we are concerned on the blow-up problem of the  Navier-Stokes equations on $\Bbb R^N$, $N\geq 2$.
$$
(NS)\left\{\aligned
&\partial_t v +(v\cdot \nabla )v=-\nabla \pi +\Delta v,\\
&\textrm{div }\, v =0,\\
&v(x,0)=v_0(x).
\endaligned
\right.
$$
Given $v_0 \in L^p_\s (\Bbb R^N)$, $p\geq N$, the local existence result due to Kato (\cite{kat2}) says that
there exists $T>0$ such that a unique solution $v\in C([0, T); L^p _\s(\Bbb R^N ))$ exists, which is smooth for $t\in (0, T)$. The global in time regularity question for the Navier-Stokes equations is a well-known millennium problem in mathematics since the pioneering paper due to Leray(\cite{ler}).
Here the suitable quantities, which corresponds to $\a_k (t)$ in the previous section, are the followings.
\bqn
\gamma_p (t)&:=&\frac{\int_{\Bbb R^N}\pi  \mathrm{div} (v |v|^{p-2})  dx}{\|v\|_{L^p} ^p},\\
\delta_p (t)&:=&\frac{\int_{\Bbb R^N } |\nabla v |^2 |v|^{p-2} dx +(p-2)\int_{\Bbb R^N} |\nabla  |v| |^2 |v|^{p-2}dx}{\|v\|_{L^p} ^p},
\eqn
and
$$
\lambda_p (t):= \gamma_p (t)-\delta_p (t)
$$
for each  $p\in [1, \infty)$, where $v(\cdot, t)\neq 0$. In the case $v(\cdot, t)=0$ we set $\gamma_p (t)=\delta_p (t)=\lambda_p (t)=0$.
The following theorem, which corresponds  to the Navier-Stokes' version of  Theorem 1.1(i), is due to Leray(pp.227, \cite{ler})(see also \cite{gig2}).
\begin{thm}
Let $p\in (N, \infty)$, and  $v\in C([0, T_*); L^p _\s(\Bbb R^N))$  be a local
smooth solution of the Navier-Stokes equations with $v_0\in L^p_\s (\Bbb R^N), v_0 \neq 0$, which blows up at $T_*$. Then, necessarily there exists a positive constant $C=C(p,N)$ such that
\bb
 \lim\inf_{t\to T_*}(T_*-t)^{\frac{p-N}{2p}}\|v(t)\|_{L^p} \geq C.
 \ee
 \end{thm}
 We state below the results corresponding to (ii) and (iii) of Theorem 1.1.
\begin{thm}
Let $p\in (N, \infty)$, and  $v\in C([0, T_*); L^p _\s(\Bbb R^N))$  be a local
smooth solution of the Navier-Stokes equations with $v_0\in L^p_\s (\Bbb R^N), v_0 \neq 0$. Then, the followings are the necessary  conditions that $T_*$ is a blow-up time of the solution $v$.
\begin{itemize}
\item[(i)]  There exists a constant $C=C(p,N)$ such that
\bb\label{1.16}
\lim\inf_{t\to T_*} \int_{0}^{t} \left\{\lambda_p (\tau)-\left[\frac{p-N}{2p}\right]\frac{1}{T_* -\tau}\right\}d\tau >-\infty.
\ee
\item[(ii)] For all $\vare_0 >1$ there exists a sequence $\{t_n \}$ with $t_n \nearrow T_*$ such that
\bb\label{1.17}
\lambda_p (t_n)\geq \left[\frac{p-N}{2p}\right]\frac{1}{T_* -t_n}-\left[\frac{(p-N)\vare_0 }{2p}\right]\frac{1 }{T_* -t_n} \left[\log \left(\frac{1}{T_* -t_n}\right)\right]^{-1}
\ee
for all $n=1,2,\cdots.$
\end{itemize}
\end{thm}
Similarly to Theorem 1.2  we establish following theorem for the Navier-Stokes equations.
\begin{thm} Let $p\in [N, \infty)$, and  $v\in C([0, T_*); L^p _\s(\Bbb R^N))$  be a local
smooth solution of the Navier-Stokes equations with $v_0\in L^p_\s (\Bbb R^N), v_0 \neq 0$..
Then, necessarily  at least one of the following three statements hold true, where the case (ii) is excluded for $p=N$.
\begin{itemize}
\item[(i)] There exists $\{ t_n \}_{n=0}^\infty$ with $t_n \nearrow T_*$  such that
\bb\label{1.18}
\lambda_p (t_n)=\left[\frac{p-N}{2p}\right]\frac{1}{T_*-t_n}\qquad \forall n=1,2,\cdots .
\ee
\item[(ii)] There exists $t_0 \in (0, T_*)$ such that
\bb\label{1.19}
 \left|\lambda_p (t)-\left[ \frac{p-N}{2p}\right] \frac{1}{T_*-t}\right| >0 \qquad \forall t\in (t_0, T_*),
\ee
and
\bb\label{1.20}
 \lim_{t\to T_*}(T_* -t)^{\frac{p-N}{2p}} \|v(t)\|_{L^p} + \int_{0}^{T_*} \left|\lambda_p (t)-\left[\frac{p-N}{2p}\right]\frac{1}{T_* -t}\right|dt <\infty .
  \ee
  \item[(iii)]
 There exists $t_0 \in (0, T_*)$ such that
\bb\label{1.21}
\lambda_p (t)> \left[ \frac{p-N}{2p}\right] \frac{1}{T_*-t}\qquad \forall t\in (t_0, T_*),
\ee
and
  \bb\label{1.22}
 \lim_{t\to T_*}(T_* -t)^{\frac{p-N}{2p}} \|v(t)\|_{L^p} =\int_{0}^{T_*} \left\{\lambda_p (t)-\left[\frac{p-N}{2p}\right]\frac{1}{T_* -t}\right\}dt =\infty,
  \ee
  and for all continuous, positive function $g$ satisfying the Osgood condition, we have
 \bb\label{1.23}
 \int_{t_1} ^{T_*} \frac{\left\{ \lambda_p  (t)-\left[\frac{p-N}{2p}\right]\frac{1}{T_* -t}\right\}}{g\left(\log\{(T_* -t)^{\frac{p-N}{2p}} \|v(t)\|_{L^p}\}\right)} dt <\infty
\ee
for $t_1\in (0, T_*)$ sufficiently close to $T_*$.
\end{itemize}
\end{thm}
{\it Remark 1.3 } Similarly to Remark 1.1 the integrability in (\ref{1.20}) implies  that there exists a sequence $\{ t_n\}$ with $t_n \nearrow T_*$ such that
\bb\label{1.24a}
\lambda_{p} (t_n)=\left[\frac{p-N}{2p}\right]\frac{1}{T_* -t_n} + o \left(\frac{1}{T_* -t_n }\left[\log \left(\frac{1}{T_* -t_n }\right)\right]^{-1}\right)  \quad\mbox{as $n \to \infty$},
\ee
which is a special case of (\ref{1.17}). We can thus further narrow down the possibilities for the behavior of $\lambda_p (t)$ near the possible blow-up time $T_*$ as follows:\\ Either (\ref{1.24a}) holds for a sequence $t_n \nearrow T_*$, or (\ref{1.21}) holds
for some $t_0<T_*$. This is  more specified than Theorem 1.4(ii).\\
\ \\
\noindent{\it Remark 1.4 } Since (ii) is excluded for $p=N$, we can summarize the main conclusion of Theorem 1.4 as follows:
If $T_*$ is the blow-up time for a local smooth solution $v\in C([0, T_*);L^N _\s(\Bbb R^N))$, then
either there exists a sequence $\{ t_n \}_{n=0}^\infty$ with $t_n \nearrow T_*$ such that
\bb\label{1.25}
\lambda_N (t_n)=0\qquad \forall n=1,2,\cdots ,
\ee
or
\bb\label{1.27}
\lim\sup_{t\to T_*} \|v(t)\|_{L^N}=\int_0 ^{T_*}| \lambda_N (t)|dt=\infty,\,\,\mbox{and}\,\, \int_{t_1} ^{T_*} \frac{|\lambda_N  (t)|}{g\left(\log( \|v(t)\|_{L^N})\right)} dt <\infty
\ee
for all positive, continuous function $g$ satisfying the Osgood condition, where $t_1$ is sufficiently close to $T_*$.

\subsection{The surface quasi-geostrophic equations}
We are concerned on the (inviscid) 2D quasi-geostrophic
 equation,
 $$
 (QG)\left\{ \aligned
 &\frac{\partial \th}{\partial t} +(v\cdot \nabla )\th =
 0,\\
 &v(x,t) =-\nabla^\bot \int_{\Bbb R^2}
 \frac{\th (y,t)}{|x-y|} dy=R^\bot \theta,\\
 &\th(x,0)=\th_0 (x),
 \endaligned
 \right.
 $$
 where $\theta (x,t)$ is a scalar function representing the temperature,
 $v(x,t)$ is the velocity field of the fluid,
 and
 $$\nabla ^\bot =(-\partial _{x_2}, \partial_{x_1} ),\quad R^\bot f =(-R_2 f, R_1 f)
 $$
with $R_j$, $j=1,2$,  are the Riesz transforms in $\Bbb R^2$ defined by
 $$ (R_j f)(x) =C\int_{\Bbb R^2} \frac{(x_j-y_j)f(y)}{|x-y|^3} dy,
 $$
where $C$ is an absolute  constant.
The system (QG) is of intensive interests recently(see e.g.
\cite{con1, con4, cor1, cor2, cor3,cor4, wu1,wu2}, and references therein),
since the equation has very similar structure to the 3D Euler
equations, and also it has direct connections to the physical
phenomena in the atmospheric science.
Since we have $L^p$ norm conservation of $\theta$,
$$
\|\theta (t)\|_{L^p}=\|\theta_0 \|_{L^p}, \qquad 0<p\leq \infty,
$$
and the local well-posedness in $W^{k,p} (\Bbb R^2)$ is known(although it is not written explicitly in the literature, the proof is straightforward, following the argument for the Euler equations due to, say, Kato\cite{kat2,kat3}).
The quantity corresponding to $\a(x,t)$ in the 3D Euler equations(\cite{con1}) is
$$
 \hat{\a}(x,t)= \sum_{j,k=1}^2\xi_j(x,t)S_{jk}  (x,t)\xi_k (x,t).
 $$
 where $\xi (x,t)=\nabla ^\bot \theta(x,t) /|\nabla ^\bot \theta (x,t)|$ is the unit tangent vector to the level, $\theta(x,t)=Const.$, and $S=(S_{jk})$ is the deformation tensor. In the case $\nabla ^\bot \theta (x,t)=0$ we set  $\hat{\a}(x,t)=0$. The natural quantity, which corresponds to $\a_k (t), \lambda_p(t)$ in the previous sections, is
\bb
\a_{k,p}(t) =\frac{-\int_{\Bbb R^2} D^k \left[(v(x,t)\cdot \nabla)\theta(x,t) \right] \cdot  D^k \theta(x,t) |D^k \theta (x,t)|^{p-2} dx}{\|D^k\theta (t)\|_{L^p}^p}.
\ee
\begin{thm}
Let  $k> 2/p +1$, and $\theta\in C([0, T_*); W^{k,p} (\Bbb R^2))$ be a solution to (QG) with $\theta_0 \in W^{k,p} (\Bbb R^2)$, $\theta_0 \neq 0$.  Then, the followings are the necessary and sufficient conditions
that $T_*$ is a blow-up time of the solution.
\begin{itemize}
\item[(i)] There exists an absolute constant $K=K(k,p)$ such that
\bb
\lim\inf_{t\to T_*} (T-t)\|D^k \theta(t)\|_{L^p}^{\frac{p+2}{kp}} \geq \frac{K}{\|\theta_0\|_{L^p}^{1-\frac{p+2}{kp}}}.
\ee
\item[(ii)]
\bb
\lim\inf_{t\to T_*} \int_{0}^{t} \left\{\a_{k,p} (\tau)-\left[\frac{kp}{p+2}\right]\frac{1}{T_* -\tau}\right\}d\tau> -\infty.
\ee
\item[(iii)]
For all $\vare_0 >1$ there exists a sequence $\{t_n \}$ with $t_n \nearrow T_*$ such that
\bb
\a_{k,p}(t_n)\geq \left[\frac{kp}{p+2}\right]\frac{1}{T_* -t_n}-\left[\frac{kp\vare_0 }{p+2}\right]\frac{1 }{T_* -t_n} \left[\log \left(\frac{1}{T_* -t_n}\right)\right]^{-1}
\ee
for all $n=1,2,\cdots.$
\item[(iv)]
For all $\vare_0 >1$ there exists a sequence $\{t_n \}$ with $t_n \nearrow T_*$ such that
\bb
\|\hat{\a}(t_n)\|_{L^\infty} \geq \frac{1}{T_* -t_n}-\frac{\vare_0  }{T_* -t_n} \left[\log \left(\frac{1}{T_* -t_n}\right)\right]^{-1}
\ee
for all $n=1,2,\cdots.$
\end{itemize}
\end{thm}

\begin{thm} Let  $k> 2/p +1$, and $\theta\in C([0, T_*); W^{k,p} (\Bbb R^2))$ be a solution to (QG) with $\theta_0 \in W^{k,p} (\Bbb R^2)$, $\theta_0 \neq 0$.  Then,  at least one of the following three statements hold true.
\begin{itemize}
\item[(i)]  There exists a sequence $\{ t_n \}_{1=0}^\infty$ with $t_n \nearrow T_*$ such that
\bb
\a _{k,p} (t_n)=\left[\frac{kp}{p+2}\right]\frac{1}{T_* -t_n} \qquad \forall n=1,2,\cdots .
\ee

\item[(ii)] There exists $t_0 \in (0, T_*)$ such that
\bb \left|\a_{k,p} (t)-\left[\frac{kp}{p+2}\right] \frac{1}{T_* -t} \right|> 0\quad \forall t\in (t_0 , T_*),
\ee
and
\bb
\lim_{t\to T_*} (T_*-t)\|D^k \theta(t) \|_{L^p}^{\frac{p+2}{kp}}+\int_{t_0}^{T_*} \left|\a_{k,p}(t)-\left[\frac{kp}{p+2}\right]\frac{1}{T_*-t}\right|dt <\infty
  \ee
  \item[(iii)] There exists $t_0 \in (0, T_*)$ such that
\bb \a_{k,p} (t) >\left[\frac{kp}{p+2}\right] \frac{1}{T_* -t} \quad \forall t\in (t_0 , T_*),
\ee
and
 \bb
\lim_{t\to T_*}(T_*-t)\|D^k \theta(t) \|_{L^p}^{\frac{p+2}{kp}}=\int_{0}^{T_*} \left\{\a_{k,p}(t)-\left[\frac{kp}{p+2}\right]\frac{1}{T_*-t}\right\}dt=\infty.
\ee
Furthermore, for any continuous, positive function $g$ satisfying the Osgood condition, we have
 \bb
 \int_{t_1} ^{T_*} \frac{\left\{ \a_{k,p}(t)-\left[\frac{kp}{p+2}\right]\frac{1}{T_* -t}\right\}}{g\left(\log\{(T_*-t)\|D^k v(t) \|_{L^p}^{\frac{p+2}{kp}} \}\right)} dt <\infty
\ee
for  $t_1\in (0, T_*)$ sufficiently close to $T_*$.
\end{itemize}
\end{thm}
{\it Remark 1.5} Similarly to Remark 1.2 the possibility of  behaviors of $\a_{k,p}(t)$ near the blow-up time $T_*$ can be summarized as follows:
Either there exists a sequence $\{ t_n\}$ with $t_n \nearrow T_*$  such that
$$
\a_{k,p} (t_n )=\left[\frac{kp}{p+2}\right]\frac{1}{T_* -t} +o\left(\frac{1}{T_* -t}\left[ \log\left(\frac{1}{T_*-t_n }\right) \right]^{-1}\right)\qquad \mbox{as}\quad n\to \infty,
$$
or there exists $t_0 \in (0, T_*)$ such that
$$
\a_{k,p} (t)> \left[\frac{kp}{p+2}\right]\frac{1}{T_* -t}\qquad\forall t\in [t_0, T_*).
$$

\section{Proof of the Main Theorems}

\subsection{The Euler equations}
We recall here the commutator estimate(\cite{kat3, kla}),
\bb\label{2.37}
\|D^k (fg)-fD^k g\|_{L^p}\leq C(\|\nabla f \|_{L^\infty} \|D^{k-1} g\|_{L^p}+\|D^k f\|_{L^p} \|g\|_{L^\infty} ),
\ee
where $1 <p<\infty$, and $k\in \Bbb N$,
and the Gagliardo-Nirenberg inequality,
\bb\label{2.38}
\|\nabla f\|_{L^\infty}\leq C \|D^k f\|_{L^p} ^{\frac{p+N}{kp}}\|f\|_{L^p}^{1-\frac{p+N}{kp}},
\ee
where $1\leq p\leq \infty$, and $k>N/p+1$, in both of which the constants $C$ depends only on $p,k$ and $N$.\\
\ \\
\noindent{\bf Proof  of Theorem 1.1 } \\
\noindent{\it \underline{Proof of part (i)}: } Given $k>N/2 +1$, we operate $D^k$ on  the first equation of (E), and then taking $L^2$ inner product it by $D^k v$.  Then we obtain after integration by part
\bq\label{2.39}
\frac12\frac{d}{dt} \|D^k v\|_{L^2}^2 &=&-\int_{\Bbb R^N}  D^k \{ (v \cdot \nabla )v\}\cdot D^k v dx\n \\
 &=&-\int_{\Bbb R^N} \left[ D^k \{ (v \cdot \nabla )v\} -(v\cdot \nabla )D^k v \right]\cdot D^k v dx\n\\
 &\leq& C\|D^k \{(v \cdot \nabla )v\} -(v\cdot \nabla)D^k v\|_{L^2} \|D^k v\|_{L^2}
 \leq C \|\nabla v\|_{L^\infty}\|D^k v\|_{L^2}^2\n \\
 &\leq& C \|D^k v (t)\|_{L^2}^{2+\frac{N+2}{2k}} \|v(t)\|_{L^2}^{1-\frac{N+2}{2k}}
 = C_{k,N} \|D^k v (t)\|_{L^2}^{2+\frac{N+2}{2k}}\|v_0\|_{L^2}^{1-\frac{N+2}{2k}}\n \\
\eq
 for an absolute constant $C_{k,N}$, in which we used the fact,
 $$
 \int_{\Bbb R^N} (v \cdot \nabla )D^k v\cdot D^k v\,dy =\frac12 \int_{\Bbb R^N} (v\cdot \nabla )|D^k v |^2dy=-\frac12  \int_{\Bbb R^N}(\mathrm{ div}\, v)|D^k v|^2 dx=0,
 $$
 and  the inequalities (\ref{2.37}) and (\ref{2.38}).
 Let us  define
 $$
X(t):= \|D^k v (t)\|_{L^2}^{\frac{N+2}{2k}}\|v_0\|_{L^2}^{1-\frac{N+2}{2k}}, \quad X_0
=\|D^k v _0\|_{L^2}^{\frac{N+2}{2k}}\|v_0\|_{L^2}^{1-\frac{N+2}{2k}},
 $$
 and set
\bb
\label{2.40}
K :=\frac{k}{(N+2)C_{k,N}}.
\ee
 Then, from (\ref{2.39}), we  compute
 $$
 \frac{d X(t)}{dt} \leq \frac{(N+2)C_{k,N}}{2k} X(t)^2=\frac{ 1}{2K} X(t)^2,
$$
which can be solved to provides us with
\bb\label{2.41}
X(t)\leq \frac{2K X_0}{2K-  X_0 t}\qquad \forall t\in [0, T_*).
\ee
Let us suppose the reverse inequality of (\ref{1.2}) holds true with constant $K$ defined by (\ref{2.40}), namely
$$
\lim\inf_{t\to T_*} (T_*-t)\|D^k v\|_{L^2} ^{\frac{N+2}{2k}}< \frac{K}{\|v_0\|_{L^2} ^{1-\frac{N+2}{2k}}}.
$$  Then there exists $t_0\in (0, T_*)$ such that
$$
(T_*-t_0 )X(t_0 ) <K.
$$
Translating the origin of  time  into $t_0$ in (\ref{2.41}), we have that for all  $t_1\in (t_0, T_*)$
\bq\label{2.42}
X(t_1)&\leq &\frac{2K X(t_0)}{2K-  (t_1 -t_0)X(t_0)}\n \\
&\leq& \frac{2K X(t_0)}{2K-  (T_* -t_0)X(t_0)}< 2X(t_0).
\eq
Passing $t_1 \nearrow T_*$ in (\ref{2.42}), we find that
\bqn
X(T_* ) &=&\|D^k v (T_*)\|_{L^2}^{\frac{N+2}{2k}}\|v_0\|_{L^2}^{1-\frac{N+2}{2k}}\\
 &\leq &2 \|D^k v(t_0)\|_{L^2}^{\frac{N+2}{2k}}\|v(t_0)\|_{L^2}^{1-\frac{N+2}{2k}}=2X(t_0)<\infty,
\eqn
 which shows that $T_*$ is not a blow-up time.  We have shown that (\ref{1.2}) is a necessary condition that $T_*$ is a blow-up time. The proof of sufficiency part is rather immediate,
  since the regularity of $v$ at $T_*$ implies  $v\in C([0, T_*];H^k_\s (\Bbb R^N))$, and hence
 $$
 \lim\inf_{t\to T_*}(T_* -t) \|D^k v(t)\|_{L^2} ^{\frac{N+2}{2k}} = \lim_{t\to T_*}(T_*-t)\|D^k v(T_*)\|_{L^2} ^{\frac{N+2}{2k}}=0,
 $$
 which is in contradiction to (\ref{1.2}).\\
 \ \\
 \noindent{\it \underline{Proof of part (ii)}:}  We write the first line of (\ref{2.39}) in the form,
\bb\label{2.43}
 \frac12\frac{d}{dt} \|D^k v\|_{L^2}^2 =-\int_{\Bbb R^N}  D^k \{ (v \cdot \nabla )v\}\cdot D^k v\, dx=\a_k (t) \|D^k v\|_{L^2}^2,
\ee
and  set
$$
Y(t):=(T_*-t )X(t)=(T_*-t) \|D^k v (t)\|_{L^2}^{\frac{N+2}{2k}}\|v_0\|_{L^2}^{1-\frac{N+2}{2k}}.
$$
Then, we compute, using (\ref{2.43}),
\bqn
\frac{d}{dt} Y(t)&=&-X(t) +\frac{N+2}{2k}\a_k (t) (T_*-t)X(t)\\
&=&\frac{N+2}{2k}\left\{ \a_k (t)- \left[\frac{2k}{N+2}\right] \frac{1}{T_* -t}\right\} Y(t).
\eqn
Integrating this over $[0,t]$, we find that
$$
 Y(t)=Y(0)\exp\left(\frac{N+2}{2k}\int_0 ^t\left\{ \a_k (\tau)- \left[\frac{2k}{N+2}\right] \frac{1}{T_* -\tau}\right\}d\tau\right),
$$
 and hence
\bq\label{2.43a}
\lefteqn{(T_*-t) \|D^k v (t)\|_{L^2}^{\frac{N+2}{2k}}=T_* \|D^k v_0\|_{L^2}^{\frac{N+2}{2k}}\times}\hspace{.5in}\n \\
&&\times\exp\left(\frac{N+2}{2k}\int_0 ^t \left\{ \a_k (\tau)- \left[\frac{2k}{N+2}\right] \frac{1}{T_* -\tau}\right\}d\tau\right).
\eq
 Thus, if
$$
\lim\inf_{t\to T_*} \int_{0}^{t} \left\{\a_k (\tau)-\left[\frac{2k}{N+2}\right]\frac{1}{T_* -\tau}\right\}d\tau =-\infty,
$$
then
 $$
 \lim\inf_{t\to T_*}(T_*-t) \|D^k v (t)\|_{L^2}^{\frac{N+2}{2k}}=0,
 $$
 and applying the result of (i), we find  that $T_*$ is not the blow-up time, and this shows that (\ref{1.3}) is a necessary condition that $T_*$ is a blow-up time.  In order to show the sufficiency  we suppose
 $v\in C([0, T_*]; H^k_\s (\Bbb R^N))$.
 Then, by (\ref{1.5})  and (\ref{2.38}) together with the energy conservation we have
\bqn
 \lefteqn{\lim\inf_{t\to T_*} \int_0 ^t \a _k (\tau )d\tau \leq \int_0 ^{T_*} |\a _k (\tau )|d\tau \leq C \int_0 ^{T_*}\|\nabla v(\tau)\|_{L^\infty} d\tau }\hspace{.5in}\n \\
  &&\leq C T_*\sup_{0<t<T_*}\|D^k v (t)\|_{L^2}^{\frac{N+2}{2k}}\|v_0\|_{L^2}^{1-\frac{N+2}{2k}}<\infty .
\eqn
 Hence,
\bqn
\lefteqn{\lim\inf_{t\to T_*} \int_{0}^{t} \left\{\a_k (\tau)-\left[\frac{2k}{N+2}\right]\frac{1}{T_* -\tau}\right\}d\tau }\hspace{.2in}\n \\
&&\leq\lim\inf_{t\to T_*} \int_{0}^{t} \a_k (\tau)d\tau-\frac{2k}{N+2}\int_0 ^{T_*} \frac{d\tau}{T_* -\tau}=-\infty,
\eqn
which is in contradiction to (\ref{1.3}).\\
 \ \\
 \noindent{\it \underline{Proof of part (iii)}:} Suppose (iii) does not hold. Then, there exists $\vare_0 >1$ and $t_0 \in (0, T_*)$ such that
\bb
\a_k (t)<\left[\frac{2k}{N+2}\right]\frac{1}{T_* -t}-\left[\frac{2k\vare_0}{N+2}\right]\frac{1 }{T_* -t} \left[\log \left(\frac{1}{T_* -t}\right)\right]^{-1}
\ee
 holds for all $t\in (t_0, T_*)$. Hence, we estimate, from (\ref{2.38}) and  (\ref{2.43a}),
 \bqn
 \lefteqn{ \int_{t_0} ^{T_*} \|\nabla v (t)\|_{L^\infty}dt\leq C \int_{t_0} ^{T_*} \|D^k v (t)\|_{L^2}^{\frac{N+2}{2k}}\|v(t)\|_{L^2}^{1-\frac{N+2}{2k}}dt=CT_*\|D^k v_0\|_{L^2}^{\frac{N+2}{2k}}\|v_0\|_{L^2}^{1-\frac{N+2}{2k}}\times}\hspace{.0in}\n \\
 &&\qquad\times \int_{t_0} ^{T_*}
 \frac{1}{T_* -t} \exp\left( \frac{N+2}{2k}\int_{0}^{t} \left\{\a_k (\tau)-\left[\frac{2k}{N+2}\right]\frac{1}{T_* -\tau}\right\}d\tau\right) dt\\
 &&\leq CT_*\|D^k v_0\|_{L^2}^{\frac{N+2}{2k}}\|v_0\|_{L^2}^{1-\frac{N+2}{2k}}\int_{t_0} ^{T_*}\frac{1}{T_* -t}\exp\left( -\int_{0} ^{t} \frac{\vare_0}{T_*-\tau}\left[\log \left(\frac{1}{T_* -\tau}\right)\right]^{-1}d\tau\right)   dt \n \\
 &&=CT_*\|D^k v_0\|_{L^2}^{\frac{N+2}{2k}}\|v_0\|_{L^2}^{1-\frac{N+2}{2k}}\int_{t_0} ^{T_*}\frac{1}{T_* -t}\left[\log \left(\frac{T_*}{T_* -t}\right)\right]^{-\vare_0} dt<\infty.
 \eqn
 Hence, thanks to the Beale -Kato-Majda criterion we find that $v\in C([0, T_*]; H^k (\Bbb R^N ))$.
  The sufficiency is immediate, since
as shown in the proof of part (ii), assumption of $v\in C([0,T_* ];H^k _\s(\Bbb R^N))$  implies that
\bqn
 \lim\sup _{t\to T_*} (T_*-t) \a_k (t)&\leq& C\lim\sup _{t\to T_*}(T_*-t) \|\nabla v(t)\|_{L^\infty}\\
  &\leq& C\lim\sup _{t\to T_*}(T_*-t) \|D^k v (T_*)\|_{L^2}^{\frac{N+2}{2k}}\|v_0\|_{L^2}^{1-\frac{N+2}{2k}}=0,
\eqn
contradicting (\ref{1.4}). \\
\ \\
 \noindent{\it \underline{Proof of part (iv)}:} We recall the well-known equation for $|\o (x,t)|$(see e.g. \cite{maj}),
 $$
 \frac{\partial}{\partial t} |\o ((X(a,t),t)|=\a (X(a,t),t) |\o(X(a,t),t)|,\quad a \in \Bbb R^3
$$
where $\{ X(a,t)\}$ is the particle trajectory generated by the velocity field $v(x,t)$. From this we compute,
$$
\frac{\partial}{\partial t} \left\{(T_* -t)|\o (X(a,t),t)|\right\}=\left(\a (X(a,t),t) -\frac{1}{T_* -t} \right) \left\{(T_* -t)|\o (X(a,t),t)|\right\}.
$$
Hence, integrating over $[0, t]$, we obtain
$$
|\o(X(a,t),t)|=\frac{T_*}{T_* -t} |\o_0 (a)|\exp\left[\int_{0} ^{t} \left(\a (X(a,\tau),\tau) -\frac{1}{T_* -\tau} \right)d\tau\right],
$$
and
$$
\|\o(t)\|_{L^\infty}\leq \frac{T_*}{T_* -t} \|\o_0 \|_{L^\infty}\exp\left[\int_{0} ^{t} \left(\|\a (\tau)\|_{L^\infty} -\frac{1}{T_* -\tau} \right)d\tau\right].
$$
From now on, repeating the argument of the proof of part(iii) above, we derive our conclusion for (iv).
$\square$\\
\ \\
\noindent{\bf Proof of Theorem 1.2 } Suppose $T_*$ is a blow-up time for
$v\in C([0, T_*); H^k (\Bbb R^N))$.
Then, at least one of the followings holds true.
 \begin{itemize}
 \item[(a)] There exists a sequence $\{t_n\}$ with $t_n \nearrow T_*$ such that
$$
\a_{k} (t_n )= \left[\frac{2k}{N+2}\right] \frac{1}{T_*-t_n} \qquad \forall n=1,2,\cdots,
$$
\item[(b)] There exists $t_0 \in (0, T_*)$ such that
$$
\left|\a_{k} (t)-\left[\frac{2k}{N+2}\right]\frac{1}{T_*-t}\right|>0 \qquad \forall t\in [t_0 , T_*),
$$
and
$$
\int_{t_0} ^{T_*} \left|\a_{k} (t)-\left[\frac{2k}{N+2}\right]\frac{1}{T_*-t}\right|dt <\infty.
$$
\item[(c)] There exists $t_0 \in (0, T_*)$ such that
$$
\a_{k} (t)-\left[\frac{2k}{N+2}\right]\frac{1}{T_*-t}<0 \qquad \forall t\in [t_0 , T_*),
$$
and
\bb\label{blow}
 \int_{t_0} ^{T_*} \left\{ \a_{k} (\tau)-\left[\frac{2k}{N+2}\right]\frac{1}{T_*-\tau}\right\} d\tau
=-\infty.
\ee
\item[(d)] There exists $t_0 \in (0, T_*)$ such that
$$
\a_{k} (t)-\left[\frac{2k}{N+2}\right]\frac{1}{T_*-t}>0 \qquad \forall t\in [t_0 , T_*),
$$
and
$$
 \int_{t_0} ^{T_*}\left\{ \a_{k} (\tau)-\left[\frac{2k}{N+2}\right]\frac{1}{T_*-\tau}\right\} d\tau
=\infty.
$$
\end{itemize}
The case (c) is eliminated, since (\ref{blow}) implies
$$\lim_{t\to T_*} (T_*-t)\|D^k v (t)\|_{L^2}^{\frac{N+2}{2k}}=0
$$
by (\ref{2.43a}), and due to Theorem 1.1 (i) $T_*$ is not a blow-up time.
 The cases (a) and (b) correspond to (i) and (ii) of Theorem 1.2 respectively, due to the formula (\ref{2.43a}).
In the case (d) we define
\bqn
 s(t)&:=&\int_{t_0} ^t  \left\{ \a_{k} (\tau)-\left[\frac{2k}{N+2}\right]\frac{1}{T_*-\tau}\right\} d\tau +\log (T_*\|D^k v(t_0)\|_{L^2} ^{\frac{N+2}{2k}})\n \\
 &=& \log\left((T_*-t)\|D^k v (t)\|_{L^2}^{\frac{N+2}{2k}}\right).
\eqn
Then, $t\mapsto s(t)$ is a monotone increasing function on $(t_0 , T_*)$ with
$s(t) \nearrow \infty$ as $t\to T_*$. Hence, we can choose $t_1\in (t_0, T_*)$
so that $s_1:= s (t_1)>1$. Furthermore,
$$s'(t)=  \a_{k} (t)-\left[\frac{2k}{N+2}\right]\frac{1}{T_*-t}>0\qquad t\in (t_0, T_*). $$
Thus, using  (\ref{2.43}), we find that for any positive, continuous function $g$ satisfying the Osgood condition the following estimate holds
$$
 \int_{t_1}^{T_*} \frac{\left|\a_{k}(t)-\left[\frac{2k}{N+2}\right]\frac{1}{T_*-t}\right|}{g\left(\log\left((T_*-t)\|D^k v (t)\|_{L^2}^{\frac{N+2}{2k}}\right)\right)}dt =\int_{t_1} ^{T_*} \frac{s'(t)}{g(s(t))}dt=\int_{s_1} ^\infty \frac{ds}{g(s)}<\infty.
$$
Thus, we have established part (iii) of Theorem 1.2
$\square$

\subsection{The Navier-Stokes equations}

In order to prove Theorem 1.4 and Theorem 1.5 we need the following lemma.
\begin{lem} Let $v\in C([0, T);L^p _\s(\Bbb R^N ))$, $N\leq p<\infty$, be a classical solution the Navier-Stokes equations in $\Bbb R^N$. Then, for all $t\in [0, T)$ there holds
\bb\label{2.58}
\|v_0\|_{L^p} \exp\left(-\int_0 ^t |\lambda_p (\tau)|\, d\tau\right)\leq \|v(t)\|_{L^p}
\leq \|v_0\|_{L^p} \exp\left(\int_0 ^t |\lambda_p (\tau)|\,d\tau\right).
\ee
\end{lem}
{\bf Proof } We take $L^2$ inner product the first equation of (NS) with $v|v|^{p-2}$, then
\bqn
\frac1p \frac{d}{dt} \|v(t)\|_{L^p} ^p&=&  \int_{\Bbb R^N} \pi (v\cdot \nabla ) |v|^{p-2} dx
 -\int_{\Bbb R^N} |\nabla v|^2 |v|^{p-2}dx  \n \\
 &&\qquad - (p-2)\int_{\Bbb R^N} |\nabla |v||^2 |v|^{p-2}dy \n \\
 &=& (\gamma_p(t) -\delta_p (t)) \|v \|_{L^p}^p=\lambda_p (t) \|v \|_{L^p}^p,
 \eqn
where we used the computation,
$$
\int _{\Bbb R^N} (\Delta v)\cdot v|v|^{p-2} dx = -\int
_{\Bbb R^N} |\nabla v|^2 |v|^{p-2}dx - (p-2)\int
_{\Bbb R^N} |\nabla |v||^2 |v|^{p-2}dx .
$$
Hence,
\bb\label{2.59}
 \frac{d}{dt} \|v(t)\|_{L^p}=\lambda_p (t)\|v \|_{L^p},
 \ee
 and
$$
\|v(t)\|_{L^p}=\|v_0\|_{L^p} \exp\left(\int_0 ^t \lambda_p (s)\,ds \right),
$$
which provides us with  (\ref{2.58}). $\square$\\
\ \\
\noindent{\bf Proof Theorem 1.4 }  \\
\noindent{\it \underline{Proof of part (i)}:}
Let us st
$$
Y(t):=(T_* -t)^{\frac{p-N}{2p}} \|v(t)\|_{L^p},
$$
and using (\ref{2.59}), we compute
$$
\frac{d Y(t)}{dt} = \left\{\lambda_p (t)-\left[\frac{p-N}{2p}\right]\frac{1}{T_*-t}\right\}Y(t),
$$
which provides us with
$$
Y(t)=Y(0) \exp \left(\int_0 ^t \left\{\lambda_p (\tau)-\left[\frac{p-N}{2p}\right]\frac{1}{T_*-\tau}\right\}d\tau\right),
$$
and
\bq\label{2.60}
\lefteqn{(T_* -t)^{\frac{p-N}{2p}}\|v(t)\|_{L^p}=T_*^{\frac{p-N}{2p}} \|v_0\|_{L^p}\times}\hspace{1.in}\n \\ &&\times\exp\left( \int_0 ^t \left\{ \lambda_p (\tau)-\left[\frac{p-N}{2p}\right] \frac{1}{T_*-\tau}\right\}d\tau\right).\n \\
\eq
Thus if
$$
\lim\inf_{t\to T_*}\int_0 ^t \left\{ \lambda_p (\tau)-\left[\frac{p-N}{2p}\right] \frac{1}{T_*-\tau}\right\}d\tau=-\infty
$$
holds, then (\ref{2.60}) implies that
$$
\lim\inf_{t\to T_*}(T_* -t)^{\frac{p-N}{2p}}\|v(t)\|_{L^p}=0,
$$
which shows that $T_*$ is not the blow-up time due to Theorem 1.3. Thus we have proved (\ref{1.16}) is a necessary condition for $T_*$ to be a blow-up time.  In order  to show sufficiency part
we assume   $v\in C([0, T_*+\vare ];L^p_\s (\Bbb R^N))$ for some $\vare >0$. Then, by standard regularity results on the Navier-Stokes equations $v(\cdot ,t) \in C^\infty _\s(\Bbb R^N)$ for all $t\in (0, T_* +\vare )$ for some $\vare >0$. We  recall that
the pressure is represented by
$$
\pi=-(\Delta )^{-1}\mathrm{ div }\,\mathrm{div}\, v\otimes v=\sum_{j,k=1}^N R_j R_k( v_j v_k),
$$
where $R_j, j=1,\cdots, N,$ are the Riesz transforms in $\Bbb R^N$(\cite{ste}). Hence, by the Calderon-Zygmund type of inequality we estimate
\bqn
\lefteqn{\left|\int_{\Bbb R^N} \pi (v\cdot \nabla ) |v|^{p-2} dx\right| =\left|\int_{\Bbb R^N} \nabla \pi \cdot v |v|^{p-2} dx\right| }\hspace{.2in}\\
&&\leq \|\nabla \pi\|_{L^p} \|v\|_{L^p} ^{p-1}
\leq  C\|v \nabla v\|_{L^p} \|v\|_{L^p} ^{p-1}\leq C\|\nabla v\|_{L^\infty} \|v\|_{L^p}^p,
\eqn
which provides us with
$|\gamma_p (t)|\leq C \|\nabla v\|_{L^\infty}$. Since $\delta _p (t) \geq 0$, we have an estimate
$$\int_0 ^{T_*} |\lambda_p (t)|dt \leq \int_0 ^{T_*} |\gamma_p (t)|dt -\int_0 ^{T_*} \delta_p (t)dt
\leq C \int_0 ^{T_*} \|\nabla v(t)\|_{L^\infty}dt <\infty.
$$
Thus  we deduce
\bqn
\lefteqn{\lim\inf_{t\to T_*}\int_0 ^t \left\{ \lambda_p (\tau)-\left[\frac{p-N}{2p}\right] \frac{1}{T_*-\tau}\right\}d\tau}\hspace{1.in}\\
&&\leq \int_0 ^{T_*} |\lambda_p (t)|dt -\frac{p-N}{2p}\int_0 ^{T_*}\frac{dt}{T_*-t}=-\infty,
\eqn
contradicting (\ref{1.16}) for $p>N$.
We have thus shown that (\ref{1.16}) is also a sufficient condition that $T_*$ is a blow-up time in our  case $p>N$.\\
\ \\
\noindent{\it \underline{Proof of part (ii)}:}
Suppose (ii) does not hold. Then, there exists $\vare_0 >1$ and $t_0 \in (0, T_*)$ such that
\bb
\a_k (t)<\left[\frac{p-N}{2p}\right]\frac{1}{T_* -t}-\left[\frac{(p-N)\vare_0}{2p}\right]\frac{1 }{T_* -t} \left[\log \left(\frac{1}{T_* -t}\right)\right]^{-1}
\ee
 holds for all $t\in (t_0 ,T_*)$. Hence, we estimate, from  (\ref{2.60}),
 \bqn
 \lefteqn{ \int_{t_0} ^{T_*} \|v (t)\|_{L^p}^{\frac{2p}{p-N}}dt=CT_*\| v_0\|_{L^p}^{\frac{2p}{p-N}}\times}\hspace{.0in}\n \\
 &&\qquad\times \int_{t_0} ^{T_*}
 \frac{1}{T_* -t} \exp\left( \frac{2p}{p-N}\int_{0}^{t} \left\{\lambda_p (\tau)-\left[\frac{p-N}{2p}\right]\frac{1}{T_* -\tau}\right\}d\tau\right) dt\\
 &&\leq CT_*\| v_0\|_{L^p}^{\frac{2p}{p-N}}\int_{t_0} ^{T_*}\frac{1}{T_* -t}\exp\left( -\int_{0} ^{t} \frac{\vare_0}{T_*-\tau}\left[\log \left(\frac{1}{T_* -\tau}\right)\right]^{-1}d\tau\right)   dt \n \\
 &&=CT_*\| v_0\|_{L^p}^{\frac{2p}{p-N}}\int_{t_0} ^{T_*}\frac{1}{T_* -t}\left[\log \left(\frac{T_*}{T_* -t}\right)\right]^{-\vare_0} dt<\infty.
 \eqn
 Hence, applying the Serrin criterion(\cite{ser,pro,oky,lad}), we find that $v\in C([0, T_*]; L^p (\Bbb R^N ))$.
  Sufficiency is immediate, since regularity at $T_*$ implies
$\lambda_p (t)\in L^\infty([t_0, T_*])$, and  thus
$$\lim\sup_{t\to T_*} (T_* -t) \lambda_p (t) =0,
$$
violating (\ref{1.17}) in the case $p>N$ . $\square$\\
\ \\
\noindent{\bf Proof Theorem 1.5 } Suppose $T_*$ is a blow-up time for
$v\in C([0, T_*); L^p (\Bbb R^N))$.
Then, at least one of the followings holds true.
 \begin{itemize}
 \item[(a)] There exists a sequence $\{t_n\}$ with $t_n \nearrow T_*$ such that
$$
\lambda_{p} (t_n )= \left[\frac{p-N}{2p}\right] \frac{1}{T_*-t_n} \qquad \forall n=1,2,\cdots,
$$
\item[(b)] There exists $t_0 \in (0, T_*)$ such that
$$
\left|\lambda_{p}(t)-\left[\frac{p-N}{2p}\right]\frac{1}{T_*-t}\right|>0 \qquad \forall t\in [t_0 , T_*),
$$
and
$$
\int_{t_0} ^{T_*} \left|\lambda_{p}(t) -\left[\frac{p-N}{2p}\right]\frac{1}{T_*-t}\right|dt <\infty.
$$
\item[(c)] There exists $t_0 \in (0, T_*)$ such that
$$
\lambda_{p}(t)-\left[\frac{p-N}{2p}\right]\frac{1}{T_*-t}<0 \qquad \forall t\in [t_0 , T_*),
$$
and
\bb\label{blo}
\int_{t_0} ^{T_*}\left\{ \lambda_{p} (\tau)-\left[\frac{p-N}{2p}\right]\frac{1}{T_*-\tau}\right\} d\tau
=-\infty.
\ee
\item[(d)] There exists $t_0 \in (0, T_*)$ such that
$$
\lambda_{p}(t)-\left[\frac{p-N}{2p}\right]\frac{1}{T_*-t}>0 \qquad \forall t\in [t_0 , T_*),
$$
and
$$
 \int_{t_0} ^{T_*} \left\{ \lambda_{p} (\tau)-\left[\frac{p-N}{2p}\right]\frac{1}{T_*-\tau}\right\} d\tau
=\infty.
$$
\end{itemize}
The case (c) is eliminated, since (\ref{blo})  implies
$$\lim_{t\to T_*} (T_*-t)^{\frac{p-N}{2p}} \|v(t)\|_{L^p} =0
$$
by (\ref{2.60}), and $T_*$ is not a blow-up time by Theorem 1.3.
 The cases (a) and (b), combined with (\ref{2.60}), correspond to (i) and (ii) of Theorem 1.4 respectively.
In the case (d) we define
\bqn
 s(t)&:=&\int_{t_0} ^t  \left\{ \lambda_{p} (\tau)-\left[\frac{p-N}{2p}\right]\frac{1}{T_*-\tau}\right\} d\tau +\log (T_* ^{\frac{p-N}{2p}}\| v(t_0)\|_{L^p} )\n \\
 &=& \log\left((T_*-t)^{\frac{p-N}{2p}}\|v (t)\|_{L^p}\right).
\eqn
Then, $t\mapsto s(t)$ is a monotone increasing function on $(t_0 , T_*)$ with
$s(t) \nearrow \infty$ as $t\to T_*$. Hence, it is legitimate to choose $t_1\in (t_0, T_*)$
so that $s_1:= s (t_1)>1$. Moreover,
$$s'(t)=  \lambda_{p}(t)-\left[\frac{p-N}{2p}\right]\frac{1}{T_*-t}>0\qquad t\in (t_0, T_*). $$
Thus, using  (\ref{2.60}), we find that for any positive, continuous function $g$ satisfying the Osgood condition the following estimate holds
$$
 \int_{t_1}^{T_*} \frac{\left|\lambda_{p}(t)-\left[\frac{p-N}{2p}\right]\frac{1}{T_*-t}\right|}{
g\left(\log\left((T_*-t)^{\frac{p-N}{2p}}\| v (t)\|_{L^p}\right)\right)}dt=\int_{t_1} ^{T_*} \frac{s'(t)}{g(s(t))}dt=\int_{s_1} ^\infty \frac{ds}{g(s)}<\infty .
$$
Thus, we established part (iii) of Theorem 1.5.
In the case  $p =N$ we can eliminate (ii), since $\int_0 ^{T_*} |\lambda_N (t)|dt <\infty$ implies $v\in L^\infty([0, T_*]; L^N (\Bbb R^N))$ due to Lemma 2.1.  Applying  the  critical case of Serrin  criterion, proved by Escauriaza, Seregin and Sverak(\cite{esc}), we deduce that the solution $v$ is regular at $T_*$. $\square$

\subsection{The surface quasi-geostrophic equations}
\noindent{\bf Proof  of Theorem 1.6 } Since the proof is similar to that of Theorem 1.1, we will be brief here, presenting only essential estimates.
 Given $k>2/p +1$, we operate $D^k$ on evolution equation part of (QG), and then taking $L^2$ inner product of it by $D^k \th |D^k \th |^{p-2}$, we obtain after integration by part
\bq\label{2.72}
\lefteqn{\frac1p\frac{d}{dt} \|D^k \th\|_{L^p}^p =-\int_{\Bbb R^2}  D^k \{ (v \cdot \nabla )\th\}\cdot D^k \th |D^k \th |^{p-2} dx}\hspace{.3in}\n \\
 &=&-\int_{\Bbb R^2} \left[ D^k \{ (v \cdot \nabla )\th \} -(v\cdot \nabla )D^k \th \right]\cdot D^k \th |D^k \th |^{p-2} dx\n\\
 &\leq& C\|D^k \{(v \cdot \nabla )\th\} -(v\cdot \nabla)D^k \th\|_{L^p} \|D^k \th\|_{L^p}^{p-1}\n \\
 &\leq& C (\|\nabla v\|_{L^\infty}\|D^k \th \|_{L^p}   + \|\nabla \th \|_{L^\infty} \| D^k v \|_{L^p})\|D^k \th\|_{L^p}^{p-1}\n \\
 &\leq & C (\|D^k v \|_{L^p} ^{\frac{p+2}{kp}} \|v\|_{L^p} ^{1-\frac{p+2}{kp}}\|D^k \th \|_{L^p}
 +\|D^k\th \|_{L^p} ^{\frac{p+2}{kp}} \|\th\|_{L^p} ^{1-\frac{p+2}{kp}}\|D^k v \|_{L^p})\|D^k \th\|_{L^p}^{p-1}\n \\
 &\leq& C \|D^k \th (t)\|_{L^p}^{p+\frac{p+2}{kp}} \|\th (t)\|_{L^p}^{1-\frac{N+2}{2k}}
 = C_{k,p} \|D^k \th (t)\|_{L^2}^{p+\frac{p+2}{kp}}\|\th_0\|_{L^2}^{1-\frac{p+2}{kp}}\n \\
\eq
 for an absolute constant $C_{k,p}$, in which we used the fact,
 $$
 \int_{\Bbb R^2} (v \cdot \nabla )D^k \th\cdot D^k \th |D^k \th |^{p-2}\,dy =\frac1p \int_{\Bbb R^2} (v\cdot \nabla )|D^k \th |^pdy=-\frac1p  \int_{\Bbb R^2}(\mathrm{ div}\, v)|D^k \th|^p dx=0,
 $$
 and  the inequalities (\ref{2.37}) and (\ref{2.38}), the Calderon-Zygmund type of inequality,
 $$
 \|D^k v\|_{L^p}=\|D^k R^\bot \th \|_{L^p}\leq C \|D^k \th \|_{L^p}, \quad \forall p\in (1, \infty),\, k\in \{0\}\cap\Bbb N,
 $$
as well as the conservation of $L^p$ norm for $\theta$.
 We set
 $$
X(t):= \|D^k \th (t)\|_{L^p}^{\frac{p+2}{kp}}\|\th_0\|_{L^p}^{1-\frac{p+2}{kp}}, \quad X_0
=\|D^k \th _0\|_{L^2}^{\frac{p+2}{kp}}\|\th_0\|_{L^p}^{1-\frac{p+2}{kp}},
 $$
and
 $$K :=\frac{kp}{2(p+2)C_{k,p}}.
 $$
 Then, from (\ref{2.72}), we  deduce
$$
 \frac{d X(t)}{dt} \leq \frac{(p+2)C_{k,p}}{kp} X(t)^2=\frac{ 1}{2K} X(t)^2.
$$
 The remaining part of the proof is the same as that of proof of Theorem 1.1 with obvious changes, and we omit it.
 $\square$\\
\ \\
\noindent{\bf Proof of Theorem 1.7 } The proof is similar to that of Theorem 1.2. We present only essential parts of it. We  note that
from the first line of (\ref{2.72})
$$
\frac{d}{dt} \|D^k \theta(t)\|_{L^p} = \a_{k,p} \|D^k \theta(t)\|_{L^p}.
 $$
 Hence, setting
 $$
 Y(t):= (T_*-t)\|D^k \theta (t)\|_{L^p}^{\frac{p+2}{kp}},
 $$
 we  compute directly to get
$$
 \frac{d Y(t)}{dt} = \frac{(p+2)}{kp} \left(   \a_{k,p} (t) -\left[\frac{kp}{p+2}\right] \frac{1}{T_*-t} \right)Y(t),
$$
which provides us with
$$
(T_*-t)\|D^k \theta (t)\|_{L^p}^{\frac{p+2}{kp}}=T_* \|D^k \theta _0\|_{L^p}^{\frac{p+2}{kp}}
\exp\left( \frac{p+2}{kp}\int_0 ^{t} \left\{   \a_{k,p} (\tau) -\left[\frac{kp}{p+2}\right] \frac{1}{T_*-\tau} \right\} d\tau\right).
$$
The remaining parts of the proof is the repetition of that of Theorem 1.2 word by word with obvious modifications, and we omit them. $\square$\\

  \end{document}